\documentclass{amsart}
\usepackage{amscd}

\newcommand{\Sln}{\mathrm{Sln}}

\newcommand{\C}{\mathcal C}
\newcommand{\cf}{\mathrm{cf}}
\newcommand{\U}{\mathcal U}
\newcommand{\V}{\mathcal V}
\newcommand{\W}{\mathcal W}
\newcommand{\rimw}{\mathrm{rim}$-$w}
\newcommand{\w}{\omega}
\newcommand{\IN}{\mathbb N}
\newcommand{\cont}{{}\hat{}}
\newcommand{\pr}{\mathrm{pr}}

\newtheorem{theorem}{Theorem}
\newtheorem{corollary}{Corollary}
\newtheorem{proposition}{Proposition}
\newtheorem{problem}{Problem}
\newtheorem{lemma}{Lemma}
\theoremstyle{definition}
\newtheorem{remark}{Remark}

\begin{document}

\title{The Suslinian number and other cardinal invariants of continua}

\author{T.Banakh, V.V.Fedorchuk, J.Nikiel and M.Tuncali}        % Enter your name between curly braces
\date{}          % Enter your date or \today between curly braces

\begin{abstract} By the {\em Suslinian number} $\Sln(X)$ of a continuum $X$ we understand the smallest cardinal
number $\kappa$ such that $X$ contains no disjoint family $\C$ of
non-degenerate subcontinua of size $|\C|\ge\kappa$. For a compact
space $X$, $\Sln(X)$ is the smallest Suslinian number of a
continuum which contains a homeomorphic copy of $X$. Our principal
result asserts that each compact space $X$ has weight
$\le\Sln(X)^+$ and is the limit of an inverse well-ordered
spectrum of length $\le \Sln(X)^+$, consisting of compacta with
weight $\le\Sln(X)$ and monotone bonding maps. Moreover,
$w(X)\le\Sln(X)$ if no $\Sln(X)^+$-Suslin tree exists. This
implies that under the Suslin Hypothesis all Suslinian continua
are metrizable, which answers a question of \cite{DNTTT1}. On the
other hand, the negation of the Suslin Hypothesis is equivalent to
the existence of a hereditarily separable non-metrizable Suslinian
continuum. If $X$ is a continuum with $\Sln(X)<2^{\aleph_0}$, then
$X$ is 1-dimensional, has rim-weight $\le\Sln(X)$ and weight
$w(X)\ge\Sln(X)$. Our main tool is the inequality
$w(X)\le\Sln(X)\cdot w(f(X))$ holding for any light map $f:X\to
Y$.
\end{abstract}
\address{T.Banakh: Instytut Matematyki, Uniwersytet Humanistyczno-Przyrodniczy Jana Kochanowskiego w Kielcach
(Poland),\newline
Department of Mathematics, Ivan Franko Lviv National University, Lviv
(Ukraine),
\newline Nipissing University, North Bay (Canada)}
\email{tbanakh@yahoo.com}
\address{V.Fedorchuk: Faculty of Mechanics and Mathematics, Lomonosov Moscow State University, Vorobjevy Gory, 1, Moscow, Russia}
\email{vvfedorchuk@gmail.com}
\address{J.Nikiel: Instytut Matematyki i Informatyki, Uniwersytet Opolski, ul. Oleska 48, 45-052 Opole, Poland}
\email{nikiel@math.uni.opole.pl}
\address{M.Tuncali: Nipissing University, North Bay (Canada)}
\email{muratt@nipissingu.ca}


\thanks{The fourth  named author is partially supported by National Science
and Engineering Research Council of Canada grants No:141066-2004.}

\maketitle

In this paper we introduce a new cardinal invariant related to the
Suslinian property of continua. By a {\em continuum} we understand
any Hausdorff compact connected space. Following \cite{Lelek}, we
define a continuuum $X$ to be {\em Suslinian} if it contains no
uncountable family of pairwise disjoint non-degenerate
subcontinua. Suslinian continua were introduced by Lelek
\cite{Lelek}. The simplest example of a Suslinian continuum is the
usual interval $[0,1]$. On the other hand, the existence of
non-metrizable Suslinian continua is a subtle problem. The
properties of such continua were considered in \cite{DNTTT1}. It
was shown in \cite{DNTTT1} that each Suslinian continuum $X$ is
perfectly normal, rim-metrizable, and 1-dimensional.  Moreover, a
locally connected Suslinian continuum has weight $\le \w_1$.

The simplest examples of non-metrizable Suslinian continua are
Suslin lines. However this class of examples has a consistency
flavour since no Suslin line exists in some models of ZFC (for
example, in models satisfying (MA$+\neg$CH)~). It turns out that
any example of a non-metrizable locally connected Suslinian
continuum necessarily has consistency nature: the existence of
such a continuum is equivalent to the existence of a Suslin line,
see \cite{DNTTT1}. This implies that under the Suslin Hypothesis
(asserting that no Suslin line exists) each locally connected
Suslinian continuum is metrizable.

It is clear that each Suslinian continuum $X$ has countable Suslin
number $c(X)$. At this point we recall the definition of some
known topological cardinal invariants. Given a topological space
$X$ let \begin{itemize}
\item
$c(X)= \sup\{|\U|:\U$ is a disjoint family of non-empty open
subsets of $X\}$ be the {\em Suslin number} of $X$;
\item $l(X)=\min\{\kappa:$ each open cover of $X$ contains a subcover of size $\le\kappa\}$
be the {\em Lindel\"of number} of $X$;
\item $d(X)=\min\{|D|:D$ is a dense set in $X\}$ be the {\em density} of $X$;
\item $hl(X)=\sup\{l(Y):Y\subset X\}$ be the {\em hereditary Lindel\"of number} of $X$;
\item $hd(X)=\sup\{d(Y):Y\subset X\}$ be the {\em hereditary density} of $X$;
\item $w(X)=\min\{|\mathcal B|:\mathcal B$ is a base of the topology of $X\}$ be the {\em weight} of $X$;
\item $\rimw(X)=\min\{\sup_{U\in\mathcal B} w(\partial U):\mathcal B$ is a base of the topology of $X\}$
be the {\em $\mathrm{rim}$-weight} of $X$.
\end{itemize}

In the context of Suslinian continua, by analogy with the Suslin number $c(X)$ it is natural to
introduce a new cardinal invariant
\begin{itemize}
\item $\Sln(X)=\sup\{|\C|:\C$ is a disjoint family of non-degenerate subcontinua of $X\}$
\end{itemize}
defined for any continuum $X$ and called the {\em Suslinian number} of $X$.
Thus a continuum $X$ is Suslinian if and only $\Sln(X)\le\aleph_0$.

It is clear that $\Sln(X)\subset \Sln(Y)$ for any pair $X\subset Y$ of continua.
It will be convenient to extend the definition of $\Sln(X)$ to all Tychonov spaces letting
$$\Sln(X)=\min\{\Sln(Y):\mbox{$Y$ is a continuum containing $X$}\}$$
for a Tychonov space $X$.

Like many other cardinal invariants the Suslinian number is monotone.

\begin{proposition} If $X$ is a Tychonov space and $Y$ is a subspace of $X$, then $\Sln(Y)\le\Sln(X)$.
\end{proposition}

The cardinal invariant $\Sln(X)$ is not trivial since it can attain any infinite value.

\begin{proposition} $\Sln(X)=c(X)=w(X)=\kappa$ for the hedgehog
$X=\{(x_\alpha)_{\alpha<\kappa}:|\{\alpha<\kappa:x_\alpha\ne0\}|\le1\}\subset[0,1]^\kappa$ with $\kappa$ needles.
\end{proposition}

Note that the each hedgehog is {\em rim-finite} in the sense that it has a base
of the topology consisting of sets with finite boundaries.
Let us remark that a rim-finite continuum $X$ with uncountable Suslinian number must be
non-metrizable (because rim-countable metrizable continua are Suslinian, see \cite{Lelek}).

The Suslinian number can not increase under monotone maps.
We recall that a map $f:X\to Y$ is {\em monotone} if $f^{-1}(y)$ is connected for any $y\in Y$.

\begin{proposition}\label{monotone} If $X$ and $Y$ are compact spaces and
$f:X\to Y$ is a surjective monotone map, then $\Sln(Y)\le\Sln(X)$.
\end{proposition}

\begin{proof} Embed $X$ in a continuum $Z$ with $\Sln(Z)=\Sln(X)$.
Consider the following equivalence relation on $Z$: $x\sim y$ if
either $x=y$ or $x,y\in X$ and $f(x)=f(y)$. Let $T=Z/_\sim$ be the
quotient space and $q:Z\to T$ be the quotient map. Since all the
equivalence classes are connected, the quotient map $q$ is
monotone. Since the preimage of a connected set under a monotone
map is connected, $\Sln(T)\le\Sln(Z)$. It remains to observe that
$Y$ can be identified with a subspace of $T$, wich yields
$\Sln(Y)\le\Sln(T)\le\Sln(Z)=\Sln(X)$.
\end{proof}

In the subsequent proof we shall refer to properties of the hyperspace $\exp(X)$ of a given compact Hausdorff space $X$. The {\em hyperspace} $\exp(X)$ of $X$ is the space of all non-empty closed subsets of $X$, endowed with the Vietoris topology. It is well known that $\exp(X)$ is a compact Hausdorff space with $w(\exp(X))=w(X)$. By $\exp_c(X)$ we denote the subspace of $\exp(X)$ consisting of subcontinua of $X$. It is easy to see that $\exp_c(X)$ is a closed subspace in $\exp(X)$.

Recall that a map $f:X\to Y$ between compact Hausdorff spaces is
called {\em light} if $f^{-1}(y)$ is zero-dimensional for each
$y\in Y$.

\begin{theorem}\label{light} If $X$ and $Y$ are compact spaces and $f:X\to Y$ is a light map,
then $w(X)\le w(Y)\cdot \Sln(X)$.
\end{theorem}

\begin{proof} Let $Z$ be a continuum such that $Z\supset X$ and $\Sln(X)=\Sln(Z)$.
Embed $Y$ into the Tychonov cube $[0,1]^{\kappa}$ where
$\kappa=w(Y)$. It follows from the Tietze-Urysohn Theorem that the
map $f$ can be extended to a map $\bar f :Z\to[0,1]^{\kappa}$.
Observe that each non-empty open set $U\subset Z$ has no one-point
component (otherwise this one-point component would be a
quasi-component and consequently $Z$ would contain a non-trivial
clopen subset which contradicts the connectedness of $Z$). Each
component of $U$ contains a non-trivial subcontinuum and
consequently, $U$ has at most $\Sln(X)$ components. Denote by
$\C_U$ the family of closures of components of $U$.

Let  $\mathcal B$ be a base for the topology of $\bar f(Z)$ with
$|\mathcal B|\le\kappa$. Finally consider the family
$\C=\bigcup_{B\in\mathcal B}\C_{\bar f^{-1}(B)}$ of closed subsets
of $X$ having size $\le\kappa$. Because of the compactness of $X$,
the inequality $w(X)\le \kappa$ will follow as soon as we prove
that the family $\C$ separates the points of $X$ in the sense that
any two distinct points $x,y\in X$ lie in disjoint elements
$C_x,C_y$ of the family $\C$.

If $f(x)\ne f(y)$, then we can find two basic subsets
$B_x,B_y\in\mathcal B$ with disjoint closures such that $f(x)\in
B_x$ and $f(y)\in B_y$. Let $C_x$ be the component of
$f^{-1}(B_x)$, containing the point $x$ and $C_y$ be the component
of $f^{-1}(B_y)$, containing the point $y$. Then
$\overline{C}_x,\overline{C}_y$ are disjoint elements of $\C$
separating the points $x,y$.

Next, suppose that $f(x)=f(y)=z$. Assuming that the family $\C$
fails to separate the points $x,y$, we will get that for each
neigborhood $W\subset [0,1]^\kappa$ of $z$ there is a subcontinuum
$K_W\subset \overline{f^{-1}(W)}$ containing both $x$ and $y$. The
net $(K_W)_W$ has a limit point $K$ in the hyperspace $\exp(Z)$ of
$Z$ (the latter means that for any neighborhood $O(K)$ of $K$ in
$\exp(Z)$ and any neighborhood $W$ of $z$ in $[0,1]^\kappa$ there
is a neighborhood $V\subset W$ of $z$ such that $K_V\in O(K)$).
This limit point $K$ is a subcontinuum of $Z$ containing the
points $x$ and $y$ and lying in $f^{-1}(z)$ which is impossible as
the latter set is zero-dimensional.
\end{proof}

The previous theorem allows us to generalize  the classical
monotone-light Factorization Theorem \cite[13.3]{Nad} asserting
that any map $f:X\to Y$ between compact Hausdorff spaces can be
represented as the (unique) composition $l\circ \mu$ of a monotone
map $\mu:X\to Z$ and a light map $l:Z\to Y$. Applying the
preceding theorem and two propositions to the calculation of the
weight of the space $Z$, we conclude that $w(Z)\le
w(Y)\cdot\Sln(Z)\le w(Y)\cdot \Sln(X)$. In such a way we obtain
the following corollary.

\begin{corollary}\label{factor} Let $f:X\to Y$ be a map between compact
spaces and $f=l\circ\mu$ be the monotone-light decomposition of
$f$ into a monotone surjective map $\mu:X\to Z$ and a light map
$l:Z\to Y$. Then the space $Z$ has weight $w(Z)\le w(Y)\cdot
\Sln(X)$ and the non-degeneracy set $N_\mu =\{z\in
Z:|\mu^{-1}(z)|>1\}$ of $\mu$ has size $|N_\mu|\le\Sln(X)$.
\end{corollary}

Applying this corollary to the constant map, we get

\begin{corollary}\label{lightzero} Each compact Haudorff space $X$ admits a monotone
map $f:X\to Z$ onto a zero-dimensional space $Z$ of weight
$w(Z)\le\Sln(X)$. In particular, each zero-dimensional compact
space $Z$ has weight $w(Z)\le\Sln(Z)$.
\end{corollary}

As an another application of Theorem~\ref{factor} we prove that
each Suslinian continuum $X$ is hereditarily decomposable, that
is, $X$ contains no indecomposable subcontinuum (a continuum $X$
is {\em indecomposable} if $X$ cannot be written as the union of
two proper non-degenerate subcontinua of $X$).
% By \'Sierpi\'nski (??) Theorem \cite{??} each metrizable indecomposable continuum $X$ has $\Sln(X)=\mathfrak c$ (because $X$ has  countinuum many composants). Consequently, each metrizable Suslinian continuum is contains no indecomposable subcontinuum.

\begin{proposition} If $X$ is a Tychonov space with $\Sln(X)\le\aleph_0$,
then all compact zero-dimensional subspaces of $X$ are metrizable and all subcontinua of $X$ are decomposable.
\end{proposition}

\begin{proof} If $Z$ is a zero-dimensional compact subset of $X$,
then $w(Z)\le \Sln(Z)\le\Sln(X)\le\aleph_0$ by the preceding corollary.

Now take any subcontinuum $C$ of $X$. Then
$\Sln(C)\le\Sln(X)\le\aleph_0$, which means that the continuum $C$
is Suslinian. Let $f:C\to [0,1]$ be any non-constant map. By
Theorem~\ref{factor}, the map $f$ can be written as the
composition $f=l\circ \mu$ of a monotone map $\mu:X\to Z$ and a
light map $l:Z\to[0,1]$ of some continuum $Z$ with $w(Z)\le
\Sln(C)\le\aleph_0$. Thus, $Z$ is a metrizable Suslinian
continuum. Such a continuum is decomposable. Otherwise, since each
indecomposable continuum has uncountably many composants (see
\cite[Theorem 7', p. 213]{Kuratowski}), it would have
$\Sln(Z)>\aleph_0$. Consequently, we  can write $Z=A\cup B$ as the
sum of two properly smaller subcontinua $A,B\subset Z$. Their
preimages $\mu^{-1}(A)$ and $\mu^{-1}(B)$ under the monotone map
$\mu$ are proper subcontinua of $C$ whose union equals $C$. This
means that the continuum $C$ is decomposable.
\end{proof}

Next we prove that the hereditary Lindel\"of number of any space $X$ is bounded from above by the Suslinian number
of $X$. For Suslinian continua this result was proved in Theorem~1 of \cite{DNTTT1}.

\begin{theorem}\label{perfect} $hl(X)\le\Sln(X)$ for any Tychonov space $X$.
\end{theorem}

\begin{proof} Let $\kappa=\Sln(X)$ and $Z\supset X$ be a continuum with $\Sln(Z)=\kappa$.

First, we prove that each singleton $\{x_0\}$, $x_0\in X$, is the
intersection of $\kappa$ open sets of $Z$. We shall construct a
transfinite sequence $(W_\alpha)_{\alpha<\alpha_0}$ of closed
neigborhoods of $x_0$ and a transfinite sequence
$(K_\alpha)_{\alpha<\alpha_0}$ of non-degenerate subcontinua of
$Z$ such that  $K_\alpha\subset \bigcap_{\beta<\alpha}W_\beta$ and
$W_\alpha\cap K_\alpha=\emptyset$ for each $\alpha<\alpha_0$. To
start the construction we choose any subcontinuum $K_0\subset
Z\setminus\{x_0\}$ and take any closed neighborhood $W_0\subset Z$
of $x_0$ missing the set $K_0$. $W_0$ is not zero-dimensional and,
since $Z$ is a continuum, we can find a subcontinuum $K_1\subset
W_0$ not containing the point $x_0$. This way we can continue the
construction until the intersection
$\bigcap_{\beta<\alpha}W_\beta$ becomes zero-dimensional. Since
$\Sln(X)=\kappa$, the construction should stop at some ordinal
$\alpha_0$ of size $|\alpha_0|\le \kappa$. For this ordinal the
intersection $Y=\bigcap_{\alpha<\alpha_0}W_\beta$ is
zero-dimensional. The compactum $Y$, being zero-dimensional,
admits a light map onto the singleton. Applying
Theorem~\ref{light}, we get $w(Y)\le\Sln(Y)\le\Sln(Z)=\kappa$.
Consequently, we can find a family $(U_\beta)_{\beta<\kappa}$ of
neighborhoods of the point $x_0$ in $Z$ such that $Y\cap
\bigcap_{\beta<\kappa}U_\beta=\{x_0\}$. Then the family
$\{W_\alpha,U_\beta:\alpha<\alpha_0,\beta<\kappa\}$ has size
$\le\kappa$ and its intersection is $\{x_0\}$.

Now, take any subspace $A\subset X$ and let $\U$ be a cover of $A$
by open subsets of $Z$. Then $\cup\U$ is an open subset of $Z$ and
$B=Z\setminus\cup\U$ is a closed set in $Z$. Consider the quotient
space $Z/B=(Z\setminus B)\cup\{B\}$ and let $q:Z\to Z/B$ be the
quotient map. Since $\Sln(Z/B)\le\Sln(Z)=\kappa$, we may apply the
previous reasoning to find a family $\V$ of open neighborhoods of
the singleton $\{B\}\in Z/B$ with $\{B\}=\cap\V$ and
$|\V|\le\kappa$. Then $\W=\{q^{-1}(V):V\in\V\}$ is a family of
size $\le\kappa$ with $\cap \W=B$. The complement $Z\setminus W$
of each $W\in\W$ is a compact subset of $Z$ which can be covered
by a finite subcover of $\U$. Therefore, the union
$\bigcup_{W\in\W}(Z\setminus W)= A$ can be covered by $\le \kappa$
elements of the cover $\U$.
\end{proof}

According to \cite[3.12.10(l)]{En}, $w(X)\le 2^{hl(X)}$ for any
compact Hausdorff space. Hence, $w(X)\le 2^{\Sln(X)}$ for any
Tychonov space. In fact, we shall prove a stronger upper bound
$w(X)\le\Sln(X)^+$.

The {\em Generalized Suslin Hypothesis} asserts that for any
regular cardinal $\kappa$ there is no $\kappa$-Suslin tree, where
a tree is called {\em $\kappa$-Suslin}  if it has height $\kappa$
but contains no chain or antichain of length $\kappa$. We recall
that the classical Suslin Hypothesis asserts that there is no
$\aleph_1$-Suslin tree.

 Below, for a cardinal $\kappa$, by $\cf(\kappa)$ we denote the
 cofinality of $\kappa$ and by $\kappa^+$  the successor cardinal of $\kappa$. We identify cardinals with initial ordinals.

\begin{theorem}\label{main} Let $X$ be a Tychonov space. Then
$w(X)\le \Sln(X)^+$. Moreover, if no $\kappa^+$-Suslin tree exists
for $\kappa=\Sln(X)$, then $w(X)\le\Sln(X)$.
\end{theorem}

\begin{proof}
Let $\kappa=\Sln(X)$ and embed $X$ into a continuum $K$ with
$\Sln(K)=\Sln(X)$. Assuming that $\kappa^+<w(X)\le w(K)$, we can
find a continuous map $f:K\to Z$ of $K$ onto a continuum $Z$ of
weight $w(Z)=\kappa^{++}$. Moreover, we may assume that the map
$f$ is monotone. Indeed, if $f$ were not monotone, then it would
factorize as $f=l\circ \mu$ with $\mu:K\to Z_1$ monotone and
$l:Z_1\to Z$ light. Then $w(Z_1)\le
w(Z)\cdot\Sln(K)=\kappa^{++}\cdot\kappa=\kappa^{++}$. Now, let us
see that the conditions $\Sln(Z)\le\kappa$ and $w(Z)=\kappa^{++}$
lead to a contradiction.

Express $Z$ as the inverse limit of a well-ordered transfinite
spectrum $\{Z_\alpha:\alpha<\kappa^{++}\}$ consisting of continua
$Z_\alpha$ with $w(Z_\alpha)\le \kappa^+$. By $p_\alpha:Z\to
Z_\alpha$, $\alpha<\kappa^{++}$, we denote the (surjective) limit
projections of the spectrum.

Consider the family $\mathcal T=\{p_\alpha^{-1}(z):z\in Z_\alpha,
\alpha<\kappa^{++}, \dim p_\alpha^{-1}(z)>0\}$ of point-preimages
which are not zero-dimensional. Endowed with the inverse inclusion
order this family forms a tree. This tree has no chains of length
more than $\kappa$. Otherwise we would obtain a strictly
decreasing sequence of length $>\kappa$ consisting of closed
subsets of $Z$, which is impossible as $hl(Z)\le\Sln(Z)=\kappa$.

The tree also contains no antichain of length $>\kappa$ since
otherwise we would construct a disjoint family of size $>\kappa$
consisting of components of some elements of $\mathcal T$.
Consequently the tree $\mathcal T$ has height $\le \kappa^+$ and
all levels of the tree have size $\le\kappa$. This implies that
the tree $\mathcal T$ contains at most $\kappa^+$ elements. Since
$\kappa^+<\kappa^{++}=\cf(\kappa^{++})$, we can find an ordinal
$\alpha<\kappa^{++}$ such that for any point $z\in Z_\alpha$ the
preimage $p_\alpha^{-1}(z)\notin\mathcal T$ is zero-dimensional.
This means that the limit projection $p_\alpha:Z\to Z_\alpha$ is
light. Applying Theorem~\ref{light}, we get a contradiction:
$w(Z)\le w(Z_\alpha)\cdot\Sln(Z)\le\kappa^+$.
\smallskip

If no $\kappa^+$-Suslin tree exists, then the tree $\mathcal T$
constructed above is not $\kappa^+$-Suslin and thus has height
$\le \kappa$. In this case we replace the condition
$w(Z)=\kappa^{++}$ by $w(Z)=\kappa^+$ and see that the proof above
gives that $w(Z)\le\Sln(Z)$.
\end{proof}

\begin{corollary} If the Generalized Suslin Hypothesis holds, then $w(X)\le\Sln(X)$ for any Tychonov space $X$.
\end{corollary}

Applying Theorem~\ref{main} to Suslinian continua, we obtain the
answer to the second part of Problem 1 of \cite{DNTTT1}.

\begin{corollary} Under the Suslin Hypothesis all Suslinian continua are metrizable.
\end{corollary}

Theorem~\ref{main} allows us to describe the structure of compacta $X$ with $w(X)>\Sln(X)$.

\begin{theorem}\label{structure} Each compact space $X$ with $w(X)>\Sln(X)$ is the inverse
limit of a well-ordered spectrum $\{Z_\alpha,\pi_\alpha^\beta,\alpha\le\beta<\Sln(X)^+\}$
consisting of compacta of weight $w(Z_\alpha)\le\Sln(X)$ and monotone bonding maps $\pi_\alpha^\beta:Z_\beta\to Z_\alpha$.
\end{theorem}

\begin{proof} Let $\kappa=\Sln(X)$. It follows from Theorem~\ref{main} that $w(X)=\kappa^+$. Therefore,
we can write $X$ as the inverse limit of a well-ordered spectrum
$\mathcal S=\{X_\alpha,p_\alpha^\beta,\alpha\le\beta<\kappa^+\}$
consisting of compacta of weight $\le\kappa$ and surjective
bonding maps. Since $hl(X)\le\Sln(X)\le\kappa$, this spectrum is
factorizable in the sense that any continuous map $f:X\to Z$ into
a compact space $Z$ of weight $w(Z)\le\kappa$ can be written as a
composition $f=f_\alpha\circ p_\alpha$ of the limit projection
$p_\alpha:X\to X_\alpha$ and a continuous map
$f_\alpha:X_\alpha\to Z$ for some ordinal $\alpha<\kappa^+$, see
\cite[3.1.6]{FC}.

For each ordinal $\alpha<\kappa^+$ let $p_\alpha=l_\alpha\circ
\mu_\alpha$ be the (unique)  monotone-light decomposition of the
limit projection $p_\alpha:X\to X_\alpha$ into a monotone map
$\mu_\alpha:X\to Z_\alpha$ and a light map $l_\alpha:Z_\alpha\to
X_\alpha$. By Proposition~\ref{monotone},
$\Sln(Z_\alpha)\le\Sln(X)\le\kappa$ and by Theorem~\ref{light},
$w(Z_\alpha)\le w(X_\alpha)\cdot\Sln(Z_\alpha)\le\kappa$. Then,
there is an ordinal $\xi(\alpha)>\alpha$ such that the monotone
map $\mu_\alpha:X\to Z_\alpha$ factorizes through
$X_{\xi(\alpha)}$ in the sense that
$\mu_\alpha=\mu_{\alpha}^{\xi(\alpha)}\circ p_{\xi(\alpha)}$ for
some map $\mu_{\alpha}^{\xi(\alpha)}:X_{\xi(\alpha)}\to Z_\alpha$.

Thus we obtain the following commutative diagram

\begin{picture}(400,130)(-70,0)

\put(160,20){$Z_\alpha$}
\put(0,100){$X$}
\put(80,100){$X_{\xi(\alpha)}$}
\put(160,100){$X_\alpha$}

\put(15,103){\vector(1,0){60}}
\put(106,103){\vector(1,0){50}}
\put(165,35){\vector(0,1){50}}
\put(97,93){\vector(1,-1){60}}
\put(15,96){\vector(2,-1){140}}

\put(35,108){\small$p_{\xi(\alpha)}$}
\put(120,108){\small $p_\alpha^{\xi(\alpha)}$}
\put(170,58){\small $l_{\alpha}$}
\put(125,70){\small $\mu_\alpha^{\xi(\alpha)}$}
\put(80,50){\small $\mu_\alpha$}

\end{picture}

Let $A$ be a cofinal subset of ordinals $<\kappa^+$ such that
$\xi(\alpha)<\beta$ for any $\alpha<\beta$ in $A$. For any
$\alpha<\beta$ in $A$ define a bonding map
$\pi_\alpha^\beta:Z_\beta\to Z_\alpha$ letting $\pi_\alpha^\beta
=\mu_\alpha^{\xi(\alpha)}\circ p^\beta_{\xi(\alpha)}\circ
l_\beta$. We claim that the map $\pi_\alpha^\beta$ is monotone.
This follows from the monotonicity of the map
$\mu_\alpha=\mu_\alpha^{\xi(\alpha)}\circ
p^\beta_{\xi(\alpha)}\circ p_\beta=\mu^{\xi(\alpha)}_\alpha\circ
p^\beta_{\xi(\alpha)}\circ
l_\beta\circ\mu_\beta=\pi^\beta_\alpha\circ\mu_\beta$. Indeed, for
any point $y\in Z_\alpha$, the preimage
$(\pi_\alpha^\beta)^{-1}(y)=\mu_\beta(\mu^{-1}_\alpha(y))$ is
connected being the image of the connected set
$\mu^{-1}_\alpha(y)$.

It is easy to see that $\pi_\alpha^\gamma=\pi_\alpha^\beta\circ\pi_\beta^\gamma$ for any ordinals $\alpha<\beta<\gamma$ in $A$, which means that $\mathcal S'=\{Z_\alpha,\pi_\alpha^\beta:\alpha,\beta\in A\}$ is an inverse spectrum. Let $Z=\lim\mathcal S'$ be the limit of this spectrum. Observe that the monotone maps $\mu_\alpha:X\to Z_\alpha$, $\alpha\in A$, induce a surjective map $\mu:X\to Z$ while the light maps $l_\alpha:Z_\alpha\to X_\alpha$, $\alpha\in A$ induce a surjective map $l:Z\to X$. Since $l_\alpha\circ\mu_\alpha=p_\alpha$ for all $\alpha\in A$, the composition $l\circ\mu:X\to X$ is the identity map of $X$. Consequently, both $l$ and $\mu$ are homeomorphisms and thus $X$ can be identified with the limit $Z$ of the spectrum $\mathcal S'$ of length $\kappa^+$ consisting of compacta of weight $\le\kappa$ and monotone bonding maps.
\end{proof}

The following particular case of Theorems~\ref{main} and
\ref{structure} answers  the remaining part of Problem 1 from
\cite{DNTTT1}.

\begin{corollary} Each non-metrizable Suslinian continuum $X$ has weight $\aleph_1$ and is
the limit of an inverse spectrum of length $\aleph_1$ consisting of metrizable Suslinian continua and monotone bonding maps.
\end{corollary}

Compacta $X$ with small Suslinian number $\Sln(X)<\mathfrak c$ share many properties of Suslinian continua.

\begin{theorem} If $X$ is a continuum with $\Sln(X)<\mathfrak c$, then $\dim X\le 1$ and
$$\mbox{rim-}w(X)\le\Sln(X)\le hl(\exp_c(X))\le w(X)\le\Sln(X)^+.$$
\end{theorem}

\begin{proof} Let $\kappa=\Sln(X)$. To show that $\rimw(X)\le\Sln(X)$, take any point
$x\in X$ and a neighborhood $U\subset X$ of $x_0$. Let
$f:X\to[0,1]$ be any function with $f(x_0)=\{0\}$ and
$f^{-1}([0,1))\subset U$. Since $\Sln(X)<\mathfrak c$, the set
$\{y\in(0,1):\dim f^{-1}(y)>0\}$ has size $\le \Sln(X)<\mathfrak
c$. Consequently, we can find a point $y\in(0,1)$ whose preimage
$f^{-1}(y)\subset Z$ is zero-dimensional. By
Corollary~\ref{lightzero},
$w(f^{-1}(y))\le\Sln(f^{-1}(y))\le\Sln(X)=\kappa$.

 Now consider the neighborhood $V=f^{-1}([0,y))$ whose boundary $\partial V$ lies in $f^{-1}(y)$
 and thus has weight $w(\partial V)\le\kappa$ and is zero-dimensional.
 This proves the inequality $\rimw(X)\le\kappa$, and  shows that the small inductive dimension ${\rm ind}(X)\le 1$.
By \cite[7.2.7]{En}, $\dim(X)\le 1$.

It remains to prove that that $\kappa\le hl(\exp_c(X))\le w(X)\le\Sln(X)^+$.
The third inequality was proved in Theorem~\ref{main} while the second inequality follows from $hl(\exp_c(X))\le
w(\exp_c(X))\le w(\exp(X))=w(X)$. Assuming that
$hl(\exp_c(X))<\kappa=\Sln(X)$, let $\lambda=hl(\exp_c(X))$ and
find a disjoint family $\C$ of size $|\C|=\lambda^+$ consisting of
non-degenerate subcontinua of $X$. This family $\C$ can be
considered as a subset of the hyperspace $\exp_c(X)$ of
subcontinua of $X$. Identify $X$ with the set of all degenerate
subcontinua in $\exp_c(X)$. Since $hl(\exp_c(X))=\lambda$, the set
$\C$ contains a subset $\C'$ of size $|\C'|=|\C|=\lambda^+$ whose
closure in $\exp_c(X)$ misses the set $X$.

We claim that $\C'$ is not a scattered subspace of $\exp_c(X)$.
Let us recall that a topological space is {\em scattered} if each
its subspace has an isolated point. It is known (and can be easily
shown) that the size of a scattered space is equal to its
hereditary Lindel\"of number. Since
$|\C'|=\lambda^+>\lambda=hl(\exp_c(X))\ge hl(\C')$, the space
$\C'$ is not scattered and thus contains a subspace $\C''$ having
no isolated point.

Now we shall construct a subset $\{C_t\}_{t\in T}\subset\C''$
indexed by elements of the binary tree
$T=\bigcup_{n\in\IN}\{0,1\}^n$ as follows. The binary tree $T$
consists of finite binary sequences. Given two binary sequences
$t=(t_0,\dots,t_n)$, $s=(s_0,\dots,s_m)$ in $T$ we write $t\le s$
if $n\le m$ and $t_i=s_i$ for all $i\le n$.

Take any distinct elements $C_0,C_1\in\C''$ and observe that the
subcontinua $C_0,C_1$ are disjoint (because the family $\mathcal
C$ is disjoint). Hence, they have open neigborhoods
$U_0,U_1\subset X$ with disjoint closures.

Assuming that for some binary sequence $s=(s_0,\dots,s_n)$ the
subcontinuum $C_s\in\C''$ and its neigborhood $U_s\subset X$ is
constructed, consider the open subset $\U_s=\{C\in\C'':C\subset
U_s\}$ of the space $\mathcal C''$ and take any two distinct (and
hence disjoint) subcontinua $C_{s\cont 0},C_{s\cont 1}\in \U_s$.
Next, choose two open neigborhoods $U_{s\cont 0},U_{s\cont
1}\subset U_s$ of $C_{s\cont 0},C_{s\cont 1}$ with disjoint
closures. This finishes the inductive step.

Now, for any infinite binary sequence $s=(s_i)$ let $C_s$ be a
cluster point of the set $\{C_{(s|n)}:n\in\IN\}$ in $\exp(X)$,
where $s|n=(s_0,\dots,s_{n-1})$. It is easy to see that
$\{C_s:s\in\{0,1\}^\w\}$ is a disjoint family of subcontinua of
$X$, lying in the closure of the set $\C''$. Since this closure
misses the set $X$, each continuum $C_s$, $s\in\{0,1\}^\w$, is
non-degenerate. Thus,
$\kappa=\Sln(X)\ge|\{C_s:s\in\{0,1\}^\w\}|=\mathfrak c$, which is
a contradiction.
\end{proof}

\begin{problem} Is $\rimw(X)\le\Sln(X)$ for any compact Hausdorff space $X$?
\end{problem}

Let us remark that all examples of non-metrizable Suslinian continua considered in the introduction or in \cite{DNTTT1} contain a copy of a Suslin line and hence fail to be hereditarily separable. However (consistent) examples of hereditarily separable Suslinian continua can be constructed as well. For such a construction we  need the following definitions and the lemma.

We recall that a surjective map $f:X\to Y$ is {\em irreducible} if $f(Z)\ne Y$ for any proper closed subset $Z$ of $X$. This is equivalent to saying that a set $D\subset X$ is dense in $X$ provided $f(X)$ is dense in $Y$.

Following \cite[III.1.15]{Fe} we call a monotone map $f:X\to Y$ between two continua {\em atomic} if for every non-degenerate subcontinuum $Z\subset Y$ the map $f|f^{-1}(Z):f^{-1}(Z)\to Z$ is irreducible. This is equivalent to saying that $\overline{D}=f^{-1}(\overline{f(D)})$ for every subset $D\subset X$ whose image $f(D)$ is dense in some non-degenerate subcontinuum of $Y$.
 An atomic map $f:X\to Y$ will be called {\em $I$-atomic} if for every $y\in Y$ the preimage is a singleton or an arc in $X$.

The following lemma will be our basic tool in the subsequent inductive construction.

\begin{lemma}\label{fed} For any non-degenerate metrizable Suslinian continuum $Y$ and any countable set $Z\subset Y$ there are a metrizable Suslinian continuum $X$ and an $I$-atomic map $f:X\to Y$ whose non-degeneracy set $N(f)=\{y\in Y:|f^{-1}|>1\}$ equals $Z$.
\end{lemma}

\begin{proof} For every $z\in Z$ fix a decreasing neighborhood base $(O_n(z))_{n\in\w}$ at $z$ such that $\overline{O_{n+1}(z)}\subset O_n(z)$ for all $n\in\w$. Let $\{q_n:n\in\w\}$ be a countable dense set in $[0,1]$. Fix a map $h_z:Y\setminus\{z\}\to[0,1]$ such that $h_z(\partial O_n(z))=\{q_n\}$ where $\partial O_n(z)$ stands for the boundary of $O_n(z)$ in $Y$. Such a choice of the map $h_z$ guarantees that $h_z(C\setminus\{z\})=[0,1]$ for any non-degenerate subcontinuum $C\subset Y$ containing $z$.

Now consider the set $X=(Y\setminus Z)\cup (Z\times[0,1])$ and the map $f:X\to Y$ which is identity on $Y\setminus Z$ and $f(z,t)=z$ for each pair $(z,t)\in Z\times[0,1]\subset X$. For every $z\in Z$ let $r_z:X\to \{z\}\times[0,1]$ be a unique map such that
\begin{itemize}
\item $r_z(y)=(z,h_z(y))$ for every $y\in Y\setminus Z\subset X$;
\item $r_z(y,t)=(z,h_z(y))$ for every $(y,t)\in (Z\setminus\{z\})\times[0,1]\subset X$;
\item $r_z(z,t)=(z,t)$ for every $t\in[0,1]$.
\end{itemize}
Endow the space $X$ with the weakest topology making the maps
$f:X\to Y$ and $r_z:X\to \{z\}\times[0,1]$, $z\in Z$, continuous.
According to \cite[III.1.2]{Fe} the obtained space $X$ is
metrizable and compact. It is easy to check that  the map $f$ is
$I$-atomic, see also \cite[III.1.15]{Fe}.

Using the atomic property of $f$ and the Suslinian property of $Y$ it is easy to check that $X$ is Suslinian too.
\end{proof}

Now, we are ready for the construction of our example. We note
that similar constructions using atomic maps have been done
before, for instance see \cite{Mac}, \cite{Mo1} and \cite{Mo2}.

\begin{theorem}\label{hersep} Under the negation of the Suslin hypothesis there exists a
hereditarily separable non-metrizable Suslinian continuum $X$.
Moreover, each non-degenerate subcontinuum of $X$ is neither
metrizable nor locally connected.
\end{theorem}

\begin{proof} Assuming the negation of the Suslin hypothesis, fix a Suslin tree $(T,\le)$ such that each node $t\in T$ has uncountably many successors in $T$ and infinitely many immediate successors in $T$.
By $h(t)$ we denote the height of a node $t\in T$ and for a countable ordinal $\alpha$ we let $T_\alpha=\{t\in T:h(t)=\alpha\}$ stand for the $\alpha$th level of $T$. For two countable ordinals $\alpha<\beta$ let $\pr^\beta_\alpha:T_\beta\to T_\alpha$ denote the map assigning to a node $t\in T_\beta$ a unique node $t'\in T_\alpha$ with $t'<t$. We may additionally assume that the tree $T$ is continuous in the sense that for any limit countable ordinal $\alpha$ and distinct nodes $t,t'\in T_\alpha$ there is $\beta<\alpha$ such that $\pr^\alpha_\beta(t)\ne\pr^\alpha_\beta(t')$.

We shall use transfinite induction to construct a well-ordered
continuous spectrum
$\{X_\alpha,\pi^\beta_\alpha:\alpha<\beta<\w_1\}$ consisting of
metrizable Suslinian continua $X_\alpha$ and atomic bonding maps
$\pi^\beta_\alpha:X_\beta\to X_\alpha$, and a sequence
$(i_\alpha:T_\alpha\to X_\alpha)_{\alpha<\w_1}$ of injective maps
such that
\begin{enumerate}
\item for any countable ordinals $\alpha<\beta$ the diagram
$$
\begin{CD}
{T_\beta}@>{i_\beta}>> {X_\beta}\\
@V{\pr^\beta_\alpha}VV @VV{\pi^\beta_\alpha}V\\
{T_\alpha}@>{i_\alpha}>> X_\alpha
\end{CD}
$$
is commutative;
\item for every $t\in T_\alpha$ the set $i_{\alpha+1}\big((\pr_\alpha^{\alpha+1})^{-1})(t)\big)$  is dense in $(\pi^{\alpha+1}_\alpha)^{-1}(i_\alpha(t))$;
\item the short projections $\pi^{\alpha+1}_\alpha:X_{\alpha+1}\to X_\alpha$ are $I$-atomic maps with non-degeneracy set $N(\pi^{\alpha+1}_\alpha)=i_{\alpha}(T_\alpha)$.
\end{enumerate}

We start the induction with a singleton $X_0$ and the injective
map $i_0:T_0\to X_0$ assigning to the root of $T$ the only point
of $X_0$. Assume that for some countable ordinal $\alpha$ the
Suslinian continua $X_\beta$, atomic bonding maps
$\pi^\beta_\gamma:X^\beta\to X_\gamma$, and injective maps
$i_\beta:T_\beta\to X_\beta$ have been constructed for all
$\gamma\le \beta<\alpha$.

If $\alpha$ is a limit ordinal, let $X_\alpha$ be the inverse
limit of the countable spectrum
$\{X_\beta,\pi^\beta_\gamma:\gamma\le\beta<\alpha\}$ and let
$\pi^\alpha_\beta:X_\alpha\to X_\beta$ stand for the limit
projections of this spectrum. They are atomic as limits of atomic
bonding maps. For every $t\in T_\alpha$ let $i_\alpha(t)$ be the
unique point of $X_\alpha$ such that
$\pi^\alpha_\beta(i_\alpha(t))=i_\beta(\pr^\alpha_\beta(t))$ for
every $\beta<\alpha$. The continuity of the tree $T$ implies that
the obtained map $i_\alpha:T_\alpha\to X_\alpha$ is injective. The
Suslinian property of $X_\alpha$ follows from that property of the
continua $X_\beta$, $\beta<\alpha$, and the atomicity of the limit
projections $\pi^\alpha_\beta$.

If $\alpha=\beta+1$ is a successor ordinal, then we can apply Lemma~\ref{fed} to find a metrizable Suslinian continuum $X_{\alpha+1}$ and an $I$-atomic map $\pi^{\alpha+1}_\alpha:X_{\alpha+1}\to X_\alpha$ whose non-degeneracy set coincides with $i_\alpha(T_\alpha)$. Thus we satisfy the condition (3) of the inductive construction. Since for every $t\in T_\alpha$ the set $\big(\pi^{\alpha+1}_\alpha\big)^{-1}(i_\alpha(t))$ is an arc in $X_{\alpha+1}$, we can define an injective map $i_{\alpha+1}:T_{\alpha+1}\to X_{\alpha+1}$ so that $\pi^{\alpha+1}_\alpha\circ i_{\alpha+1}=i_\alpha\circ\pr^{\alpha+1}_\alpha$ and $i_{\alpha+1}$ satisfies the condition (2) of the inductive construction.

After completing the inductive construction, consider the inverse limit $X$ of the spectrum $\mathcal S=\{X_\alpha,\pi^\alpha_\beta:\beta<\alpha<\w_1\}$. Using the atomicity of the bonding projections, one can check that the limit projections $\pi_\alpha:X\to X_\alpha$ are atomic as well.

Now, we establish the desired properties of the continuum $X$.
First, we show that each non-degenerate subcontinuum $C$ of $X$ is
not metrizable and not locally connected. Let $\alpha$ be the
smallest ordinal such that $|\pi_\alpha(C)|>1$. The continuity of
the spectrum $\mathcal S$ implies that $\alpha=\beta+1$ for some
ordinal $\beta$. Then $\pi_\beta(C)$ is a singleton and hence
$\pi_\beta(C)\subset i_\beta(T_\beta)$ (otherwise $C$ would be a
singleton). Let $t\in T_\beta$ be a node of $T$ with
$\pi_\beta(C)=\{i_\beta(t)\}$. It follows that $\pi_\alpha(C)$ is
a non-degenerate subcontinuum of the arc
$A_t=(\pi_\beta^\alpha)^{-1}(i_\beta(t))$. The density of
$i_\alpha(T_\alpha)$ in $A_t$ implies the existence of a node
$t'\in T_\alpha$ with $i_\alpha(t')\in \pi_\alpha(C)$. The
atomicity of the projection $\pr_\alpha$ implies that the
continuum $C=\pi_\alpha^{-1}(\pi_\alpha(C))$ contains the
subcontinuum $\pr^{-1}_\alpha(i_\alpha(t'))$ which is not
metrizable (because $t'$ has uncountably many successors in the
tree $T$). Consequently, $C$ is not metrizable as well.

To show that $C$ is not locally connected, assume the converse
and, given any two distinct points $x,x'\in
\pr^{-1}(\alpha)(i_\alpha(t'))$, find a closed connected
neighborhood $U\subset C$ of $x$ with $x'\notin U$. Since
$\pr^{-1}_\alpha(i_\alpha(t'))$ is nowhere dense in $C$, the set
$U$ has non-degenerate projection $\pr_\alpha(U)$. Then the
atomicity of $\pr_\alpha$ implies that
$x'\in\pr^{-1}_\alpha(\pr_\alpha(U))=U$, which is a contradiction.

Next, we shall prove that the continuum $X$ is Suslinian. Take any
family  $\C$ of pairwise disjoint non-degenerate subcontinua in
$X$. Repeating the preceding argument, for every $C\in\C$ we can
find a countable  ordinal $\alpha$ and a node $t_C\in T_{\alpha}$
such that $C\supset \pi^{-1}_{\alpha}(i_\alpha(t_C))$. It follows
that the nodes $t_C$, $C\in\C$, are pairwise incomparable in $T$
(otherwise the family $\C$ would contain two intersecting
continua). Since $T$ is a Suslin tree, the antichain
$\{t_C:C\in\mathcal C\}$ is at most countable and so is the family
$\C$, witnessing the Suslinian property of $X$.

It remains to check that the continuum $X$ is hereditarily
separable. By \cite[3.12.9]{En} it suffices to prove that each
closed subspace $F$ of $X$ in separable. By Theorem~\ref{perfect},
the continuum $X$, being Suslinian, is perfectly normal and hence
$F=\pi^{-1}_\alpha(\pi_\alpha(F))$ for some countable ordinal
$\alpha$. Let $Z=\pr_\alpha(F)$. Since
$$F=\pi_{\alpha}^{-1}(Z\setminus i_\alpha(T_\alpha))\cup
\bigcup_{z\in Z\cap i_\alpha(T_\alpha)}\pi_\alpha^{-1}(z)$$ and
$\pi_{\alpha}^{-1}(Z\setminus i_\alpha(T_\alpha))$ is homeomorphic
to the metrizable separable space $Z\setminus i_\alpha(T_\alpha)$,
it remains to check that for every $z\in i_\alpha(T_\alpha)$ the
continuum $\pi^{-1}_\alpha(z)$ is separable. Consider the arc
$A=\pi^{\alpha+1}_\alpha(z)$ in $X_{\alpha+1}$ and observe that
$D=A\setminus i_{\alpha+1}(T_{\alpha+1})$ is a dense subspace of
$A$. It follows from the construction that
$\pi_{\alpha+1}^{-1}(D)$ is a topological copy of $D$, dense in
$\pi_{\alpha+1}^{-1}(A)=\pi_{\alpha}^{-1}(z)$. Therefore, the
continuum $\pi^{-1}_\alpha(z)$ is separable.
\end{proof}

We do not know if the preceding theorem can be generalized to higher cardinals.

\begin{problem} Does the existence of a $\kappa^+$-Suslin tree imply the existence of a continuum $X$ with $hd(X)\le\Sln(X)=\kappa<w(X)$?
\end{problem}

\begin{remark} The existence of a $\kappa^+$-Suslin tree is equivalent to the existence
of a linearly ordered continuum $X$ with
$\kappa=\Sln(X)=c(X)<d(X)=w(X)=\kappa^+$.
\end{remark}

The non-metrizable hereditarily separable Suslinian continuum
constructed in Theorem~\ref{hersep} is very far from being locally
connected. In \cite{DNTTT2}, it was proved that separable
homogeneous Suslinian continua are metrizable. This encourages to
remind the following question of \cite{DNTTT1}.

\begin{problem} Is each locally connected (hereditarily) separable Suslinian continuum metrizable?

\end{problem}

%\begin{problem} Does any compact Hausdorff space $X$ embed into a locally connected continuum %$Y$ with $\Sln(Y)=\Sln(X)$?
%\end{problem}

% Set the ending of a LaTeX document
\end{document}